\input amstex.tex
\input amsppt.sty
\magnification=\magstep1
\vsize=22.2truecm
\baselineskip=16truept
\nologo
\TagsOnRight

\def\Z{\Bbb Z}

\def\Q{\Bbb Q}

\def\l{\left}
\def\r{\right}
\def\bg{\bigg}
\def\({\bg(}
\def\[{\bg[}
\def\){\bg)}
\def\]{\bg]}
\def\t{\text}
\def\f{\frac}
\def\mo{\roman{mod}}

\def\eq{\equiv}
\def\cs{\ldots}
\def\ls{\leqslant}
\def\gs{\geqslant}
\def\al{\alpha}
\def\ve{\varepsilon}

\def\Proof{\noindent{\it Proof}}

\def\Ack{\medskip\noindent {\bf Acknowledgment}}
\topmatter \hbox{J. Number Theory 124(2007), no.\,1, 57--61.}
\medskip
\title Simple arguments on consecutive power residues\endtitle
\author Zhi-Wei Sun\endauthor
\affil Department of Mathematics, Nanjing University
\\Nanjing 210093, People's Republic of China
\\{\tt zwsun\@nju.edu.cn}
\\{http://math.nju.edu.cn/${}^\sim$zwsun}
\endaffil
\abstract By some extremely simple arguments, we point out the
following: (i) If $n$ is the least positive $k$th power
non-residue modulo a positive integer $m$, then the greatest
number of consecutive $k$th power residues mod $m$ is smaller than
$m/n$. (ii) Let $O_K$ be the ring of algebraic integers in a
quadratic field $K=\Q(\sqrt d)$ with $d\in\{-1,-2,-3,-7,-11\}$.
Then, for any irreducible $\pi\in O_K$ and positive integer $k$ not
relatively prime to $\pi\bar\pi-1$, there exists a $k$th power
non-residue $\omega\in O_K$ modulo $\pi$ such that
$|\omega|<\sqrt{|\pi|}+0.65$.
\endabstract
\thanks 2000 {\it Mathematics Subject Classification}. Primary 11A15;
Secondary 05A19, 11A07, 11R11.
\newline\indent
Partially supported by the National Science Fund for Distinguished
Young Scholars (No. 10425103) and a Key Program of NSF
(No. 10331020) in China.
\endthanks
\endtopmatter

\document

 For an integer $a$ relatively prime to a positive integer $m$,
if the congruence $x^k\eq a\ (\mo\ m)$ is solvable then $a$ is
said to be a $k$th power residue mod $m$, otherwise $a$ is called
a $k$th power non-residue mod $m$. The theory of power residues
(cf. [L]) plays a central role in number theory.

In this short note we aim to show that some classical topics on
power residues can be handled just by some extremely simple
observations.

Our first observation concerns the least positive $k$th power
non-residue modulo a positive integer.

\proclaim{Theorem 1} {\rm (i)} Suppose that $n=n_k(m)$ is
the least positive $k$th power non-residue modulo a positive integer
$m$. Then the greatest number $R=R_k(m)$ of consecutive $k$th
power residues mod $m$ is smaller than $m/n$, consequently $n<\sqrt
m+1/2$ if $m$ is a prime.

{\rm (ii)} Let $p$ be an odd prime, and let $k$ be a positive
integer with $\gcd(k,p-1)>1$. Provided that $-1$ is a $k$th power
residue mod $p$
$($i.e., $(p-1)/\gcd(k,p-1)$ is even$)$ and that $n_k(p)\not=2$, we have
$n_k(p)<\sqrt{p/2}+1/4.$
\endproclaim

Let $p$ be an odd prime. That $n_2(p)<\sqrt p+1$ was first pointed
out by Gauss. Using sophisticated analytic tools, A. Granville, R.
A. Mollin and H. C. Williams [GMW] proved that if $d>3705$ is a
discriminant of a quadratic number field then the Kronecker symbol
$(\f dq)$ is $-1$ for some prime $q<\sqrt d/2$; in particular, if
$p>3705$ is a prime with $p\eq1\ (\mo\ 4)$ then
$n_2(p)<\sqrt{p/4}$.

Let $k$ be a positive integer with $\gcd(k,p-1)>1$.
By a complicated
elementary method, R. H. Hudson [Hud] proved
that $n_k(p)<\sqrt{p/3}+2$ if $p\not=23,71$.
Let $N_k(p)$ denote the the greatest
number of consecutive $k$th power non-residues mod $p$.
By modifying our proof of Theorem 1 slightly, we can also show the following inequalities:
$$N_k(p)<\f{p-1}{n_k(p)-1},\ \ R_k(p)\cdot\min\{R_k(p),N_k(p)\}<p,\ \ R_2(p)N_2(p)<p.$$
A. Brauer [B] proved that $\max\{R_2(p),N_2(p)\}<\sqrt p$
for each prime $p\eq3\ (\mo\ 4)$.
By a very sophisticated elementary approach, P. Hummel [Hum] confirmed in 2003
a conjecture of I. Schur by showing that $N_2(p)<\sqrt p$ except $p=13$.

\medskip
\noindent{\it Proof of Theorem 1}. (i) Suppose that all of
$a+1,\cs,a+R$ are $k$th power residues mod $m$ where $a\in\Z$. Let
$q$ be the least integer greater than $an/m$. For any
$i\in\{1,\cs,R\}$, $(a+i)n-mq$ is a $k$th power non-residue mod
$m$ and hence
$$(a+i)n-mq\gs0\Longrightarrow (a+i)n-mq\gs n\Longrightarrow(a+i-1)n-mq\gs0.$$
As $an-mq\not\gs0$, we must have $(a+R)n-mq\not\gs0$ and thus
$nR<mq-an\ls m$. If $m$ is a prime $p$, then
 $1,\cs,n-1$ are $k$th power residues mod $p$,
therefore $n(n-1)\ls nR\ls p-1<p-1/4$ and hence $n-1/2<\sqrt p$.

(ii) Write $p=2nq+r$ with $q,r\in\Z$ and $0<|r|<n=n_k(p)$. As
$2nq=p-r$ is a $k$th power residue mod $p$,\, $q$ must be a $k$th
power non-residue mod $p$ and hence $q\gs n$ since $q>0$.
Therefore $p\gs 2n^2-(n-1)$ and thus $n-1/4<\sqrt{p/2}$. \qed
\medskip

Our second observation is the following new result established by
our simple method used in the proof of Theorem 1(ii).

\proclaim{Theorem 2} Let $K=\Q(\sqrt d)$ be a quadratic field with
$d\in\{-1,-2,-3,-7,-11\}$, and let $O_K$ be the ring of algebraic
integers in $K$. Let $\pi$ be any irreducible element of $O_K$,
and let $k$ be a positive integer with $\gcd(k,N(\pi)-1)>1$ where
$N(\pi)=\pi\bar{\pi}$ is the norm of $\pi$ with respect to the
field extension $K/\Q$. Then there is a $k$th power non-residue
$\omega\in O_K$ modulo $\pi$ with $|\omega|<\sqrt{|\pi|}+0.65$.
\endproclaim

\Proof. It is well known that $O_K$ is an Euclidean domain with
respect to the norm $N:\ O_K\to\{0,1,2,\cs\}$. (See, e.g.,
[ELS].) Thus $O_K$ is a principle ideal domain and $O_K/(\pi)$ is
a field with $N(\pi)=|\pi|^2$ elements. If $\al\in O_K$ is a $k$th
power residue mod $\pi$, then $\al^{(N(\pi)-1)/d}\eq1\ (\mo\ \pi)$
where $d=\gcd(k,N(\pi)-1)>1$. As the congruence $x^n\eq 1\ (\mo\
\pi)$ over $O_K$ has at most $n$ solutions, there are $k$th power
non-residues modulo $\pi$.

Let $\omega\in O_K$ be a $k$th power non-residue mod $\pi$ with
minimal norm. Then $N(\omega)<N(\pi)$ because
$N(\omega-\eta\pi)<N(\pi)$ for a suitable $\eta\in O_K$.

Choose $\beta,\gamma\in O_K$ so that $\pi=\beta \omega+\gamma$ and
$N(\gamma)<N(\omega)<N(\pi)$. If $\gamma=0$, then $\omega$ is a
unit since $\pi$ is irreducible, hence $N(\omega)=1<|\pi|=\sqrt
{N(\pi)}$ and so $|\omega|<\sqrt{|\pi|}$.

Now assume that $\gamma\not=0$. Then $\pi\nmid \gamma$ since
$N(\pi)\nmid N(\gamma)$. As $N(-\gamma)<N(\omega)$, $\beta\omega
=\pi-\gamma$ is a $k$th power residue mod $\pi$ and hence $\beta$
must be a $k$th power non-residue mod $\pi$. So $N(\beta)\gs
N(\omega)$, i.e., $|\beta|\gs|\omega|$. Note also that $|\gamma|<|\omega|$
since $N(\gamma)<N(\omega)$.
Therefore $|\pi|\gs|\beta|\cdot|\omega|-|\gamma|>|\omega|^2-|\omega|$. As
$$|\omega|=\sqrt{N(\omega)}\gs\sqrt2=\f{c^2}{2c-1}$$
with $c=\sqrt 2-\sqrt{2-\sqrt2}=0.6488...,$
we have
$$|\pi|>|\omega|^2-|\omega|\gs|\omega|^2-2c|\omega|+c^2=(|\omega|-c)^2$$
and hence
$|\omega|<\sqrt{|\pi|}+c<\sqrt{|\pi|}+0.65$.
This concludes the proof. \qed

\medskip

Concerning quadratic residues and non-residues, the law of quadratic reciprocity
plays a central role.
The general version of Gauss' lemma (cf. [S]),
Euler's version of the law (cf.
Proposition 5.3.5 of [IR]) and Scholz's proof of it (cf. [D, pp.\,70-73]) via
Gauss' lemma (this proof was rediscovered by the author in 2003),
lead us to give our third theorem.

\proclaim{Theorem 3} Let $a,b,m,n$ be positive integers with
$m-\ve n=2ab$ and $\gcd(a,m)=\gcd(a,n)=1$, where $\ve$ is $1$ or
$-1$. Then we have the identity
$$r_m(a)-\ve r_n(a)=\l\lfloor \f a2\r\rfloor b,$$
where
$$r_l(a)=\bg|\l\{0<r<\f l2:\ r\in\Z\ \t{and}\
\l\{\f{ar}l\r\}>\f12\r\}\bg|\quad\ \t{for}\ \ l=1,2,3,\ldots,$$
and $\lfloor\al\rfloor$ and $\{\al\}$ denote
the integral part and the fractional part of a real number $\al$ respectively.
\endproclaim

If the condition $\gcd(a,m)=\gcd(a,n)=1$ in Theorem 3 is cancelled, then
we can refine our proof of Theorem 3 to yield the following result:
$$\align&\bg|\bg\{0<r<\f m2:\ r\in\Z,\ \bg\{\f{ar}m\bg\}\gs\f12\bg\}\bg|
-\ve\bg|\bg\{0<r<\f n2:\ r\in\Z,\ \bg\{\f{ar}n\bg\}\gs\f12\bg\}\bg|
\\&\qquad\quad=\cases\lfloor a/2\rfloor b-\lfloor\gcd(a,n)/2\rfloor&\t{if}\ \ve=-1,\ \t{and}\ n/\gcd(a,n)\ \t{is odd};
\\\lfloor a/2\rfloor b&\t{otherwise}.
\endcases
\endalign$$

\noindent{\it Proof of Theorem 3}. As $\gcd(2a,m)\mid 2$, we have $2a\nmid sm$ for each
positive integer $s<a$. Clearly
 $$\align r_m(a)=&\sum_{s=0}^{a-1}\bg|\bg\{\f s{2a}m<r<\f{s+1}{2a}m:
 \ r\in\Z\ \t{and}\ \bg\{\f{ar}m\bg\}>\f12\bg\}\bg|
 \\=&\sum_{s=0}^{a-1}\bg|\bg\{r\in\Z:\ \f s2<\f{ar}m<\f{s+1}2
 \ \t{and}\ \bg\{\f{ar}m\bg\}>\f12\bg\}\bg|
 =\sum\Sb 0\ls s<a\\2\nmid s\endSb \Delta_s(m),\endalign$$
where
$$\align\Delta_s(m)=&\bg|\bg\{r\in\Z:\ \f{s}{2a}m<r<\f{s+1}{2a}m\bg\}\bg|
\\=&\bg|\bg\{r\in\Z:\ \ve\f s{2a}n<r-bs<\ve\f{s+1}{2a}n+b\bg\}\bg|
\\=&\bg|\bg\{x\in\Z:\ \ve\f s{2a}n<x<\ve\f{s+1}{2a}n+b\bg\}\bg|.
\endalign$$
Similarly, $2a\nmid sn$ for every positive integer $s<a$, and
 $$r_n(a)=\sum\Sb 0\ls s<a\\2\nmid s\endSb\Delta_s(n)$$
with
$$\Delta_s(n)=\bg|\bg\{r\in\Z:\ \f{s}{2a}n<r<\f{s+1}{2a}n\bg\}\bg|
=\bg|\bg\{x\in\Z:\ -\f{s+1}{2a}n<x<-\f{s}{2a}n\bg\}\bg|.$$
For any positive odd integer $s<a$, we have $2a\nmid (s+1)n$
(since $a\nmid s+1$ if $s<a-1$, and $2\nmid n$ if $a=s+1\eq0\ (\mo\ 2)$), hence
$$\Delta_s(m)-\ve\Delta_s(n)
=\bg|\bg\{x\in\Z:\ \ve\f{s+1}{2a}n<x<\ve\f{s+1}{2a}n+b\bg\}\bg|=b.$$
Therefore
$$r_m(a)-\ve r_n(a)=\sum\Sb 0<s<a\\2\nmid s\endSb(\Delta_s(m)-\ve\Delta_s(n))
=|\{0<s<a:\ 2\nmid s\}|\times b=\l\lfloor\f a2\r\rfloor b.$$ This
proves Theorem 3. \qed

\Ack. The author thanks the referee for his/her helpful comments.

\enddocument